\documentclass[final]{dmtcs}

\newtheorem{theorem}{Theorem}
\newtheorem{definition}{Definition}
\newtheorem{proposition}{Proposition}

\usepackage{color}
\definecolor{light-gray}{gray}{0.4}

\usepackage[latin1]{inputenc}
\usepackage{subfigure}

\usepackage{amsmath,amssymb,mathrsfs}

\newcommand{\plug}[1]{\stackrel{#1}{\blacktriangleleft}}

\author{Pawel Blasiak\thanks{Email: \email{Pawel.Blasiak@ifj.edu.pl}}}

\title[Combinatorial Route to Algebra:\\The Art of Composition \& Decomposition]{Combinatorial Route to Algebra:\\The Art of Composition \& Decomposition}

\address{H. Niewodnicza\'nski Institute of
Nuclear Physics, Polish Academy of Sciences, Krak\'ow, Poland}

\keywords{composition/decomposition of objects, multiplication/co-multiplication, combinatorics, algebra}

\received{2009-07-31}
\revised{2010-03-23}
\accepted{2010-08-10}
\begin{document}
\publicationdetails{12:2}{2010}{381}{400}
\maketitle
\begin{abstract}
We consider a general concept of composition and decomposition of objects, and discuss a few natural properties one may expect from a reasonable choice thereof. It will be demonstrated how this leads to multiplication and co-multiplication laws, thereby providing a generic scheme furnishing combinatorial classes with an algebraic structure. The paper is meant as a gentle introduction to the concepts of composition and decomposition with the emphasis on combinatorial origin of the ensuing algebraic constructions.
\end{abstract}

\section{Introduction}

A great deal of concrete examples of abstract algebraic structures are based on combinatorial constructions. 
Their advantage comes from simplicity steaming from the use of intuitive notions of enumeration, composition and decomposition which oftentimes provide insightful interpretations and neat pictorial arguments. In the present paper we are interested in clarifying the concept of composition/decomposition and the development of a general scheme which would furnish combinatorial objects with an algebraic structure.

In recent years the subject of combinatorics has grown into maturity. It gained a solid foundation in terms of combinatorial classes consisting of objects having size, subject to various constructions transforming classes one into another, see \textit{e.g.}
\cite{FlajoletBook,BergeronBook,GouldenBook}. Here, we will augment this framework by considering objects which may \emph{compose} and \emph{decompose} according to some internal law. This idea was pioneered by G.-C. Rota who considered \emph{monoidal} composition rules and introduced the concept of section coefficients showing how they lead to the co-algebra and bi-algebra structures \cite{JoniRota}. It was further given a firm foundation by A. Joyal \cite{Jo81} on the grounds of the theory of species. In the present paper we discuss this approach from a modern perspective and generalize the concept of composition in order to give a proper account for indeterminate (\emph{non-monoidal}) composition laws according to which two given objects may combine in more than one way.
Moreover, we will provide a few natural conditions one might expect from a reasonable composition/decomposition rule, and show how they lead to algebra, co-algebra, bi-algebra  and Hopf algebra structures \cite{Bo89,SweedlerBook,Ab80}. Our treatment has the virtue of a direct scheme translating combinatorial structures into algebraic ones -- it only has to be to checked whether the law of composition/decomposition of objects satisfies certain properties. We note that these ideas, however implicit in the construction of some instances of bi-algebras, have never been explicitly exposed in full generality. We will illustrate this framework on three examples of classical combinatorial structures: words, graphs and trees. For words we will provide three different composition/decomposition rules leading to the free algebra, symmetric algebra and shuffle algebra (the latter one with a non-monoidal composition law) \cite{LothaireBook,ReutenauerBook}. In the case of graphs, except of trivial rules leading to the commutative and co-commutative algebra, we will also describe the Connes-Krimer algebra of trees \cite{Kreimer1998,ConnesKreimer1998}. One can find many other examples of monoidal composition laws in the seminal paper \cite{JoniRota}; for instances of non-monoidal rules see \textit{e.g.} \cite{GrLa89,BlDuHoPeSo10}. A comprehensive survey of a recent development of the subject with an eye on combinatorial methods can be found in \cite{CartierHopfPrimer}.

The paper is written as a self-contained tutorial on the combinatorial concepts of composition and decomposition explaining how they give rise to algebraic structures. We start in Section~\ref{Preliminaries} by briefly recalling the notions of multiset and combinatorial class. In Section~\ref{Sect-Comp-Decomp} we precise the notion of composition/decomposition and discuss a choice of general conditions which lead to the construction of algebraic structures in Section~\ref{Sect-Algebra}. Finally, in Section~\ref{Sect-Examples} we illustrate this general scheme on a few concrete examples.

\section{Preliminaries}\label{Preliminaries}

\subsection{Multiset}

A basic object of our study is a \emph{multiset}. It differs from a set by allowing multiple copies of elements and formally can be defined as a pair $(A,m)$, where $A$ is a set and $m:A\longrightarrow\naturals_{\geqslant1}$ is a function counting multiplicities of elements.\footnote{For equivalent definition of a multiset based on the $\textsc{Seq}$ construction subject to appropriate equivalence relation see \cite{FlajoletBook} p. 26.} For example the multiset $\{a,a,b,c,c,c\}$ is described by the underlying set $A=\{a,b,c\}$ and the multiplicity function $m(a)=2$, $m(b)=1$ and $m(c)=3$. Note that each set is a multiset with multiplicities of all elements equal to one. It is a usual practice to drop the multiplicity function $m$ in the denotation of a multiset $(A,m)$ and simply write $A$ as its character should be evident from the context (in the following we will mainly deal with multisets!). Extension of the conventional set-theoretical operations to multisets is straightforward by taking into account copies of elements.
Accordingly, the sum of two multisets $(A,m_A)$ and $(B,m_B)$ is the multiset $A\uplus B=(A\cup B,m_{A\cup B})$ where 
\begin{eqnarray}
m_{A\cup B}(x)=\left\{\begin{array}{ll} m_A(x)+m_B(x)\,, &\ \  x\in A\cap B  \\m_A(x)\,,&\ \ x\in A-B \\m_B(x)\,,&\ \ x\in B-A  \end{array}\right.\ ,
\end{eqnarray}
whilst the product is defined as $A\times B=(A\times B,m_{A\times B})$, where
$m_{A\times B}((a,b))= m_A(a)\cdot m_B(b)$.
We note that one should be cautious when comparing multisets and not forget that equality involves the coincidence of the underlying sets as well as the multiplicities of the corresponding elements. Similarly, inclusion of multisets $(A,m_A)\subset (B,m_B)$ should be understood as the inclusion of the underlying sets $A\subset B$ with the additional condition $m_A(x)\leqslant m_B(x)$ for $x\in A$.

\subsection{Combinatorial class}

In the paper we will be concerned with the concept of \emph{combinatorial class} $\mathcal{C}$ which is a denumerable collection of objects. Usually, it is equipped with the notion of \emph{size} $|\cdot|:\mathcal{C}\longrightarrow\naturals$ which counts some characteristic carried by objects in the class, \emph{e.g.} the number of elements they are build of. The size function divides $\mathcal{C}$ into disjoint subclasses $\mathcal{C}_n=\{\varGamma\in\mathcal{C}:|\varGamma|=n\}$ composed of objects of size $n$ only. Clearly, we have $\mathcal{C}=\bigcup_{n\in\naturals}$.
A typical problem in combinatorial analysis consists in classifying objects according to the size and counting the number of elements in $\mathcal{C}_n$.

In the sequel we will often use the \emph{multiset construction}. For a given combinatorial class $\mathcal{C}$ it defines a new class $\textsc{MSet}(\mathcal{C})$ whose objects are multisets of elements taken from $\mathcal{C}$. We note that the size of $\varGamma\in\textsc{MSet}(\mathcal{C})$ is canonically defined as the sum of sizes of all its elements, \textit{i.e.} $|\varGamma|=\sum_{\gamma\in\varGamma}|\gamma|$.

\section{Combinatorial Composition \& Decomposition}\label{Sect-Comp-Decomp}

We will consider a combinatorial class $\mathcal{C}$ consisting of objects which can \emph{compose} and \emph{decompose} within the class. In this section we precisely define both concepts and discuss a few natural conditions one might expect from a reasonable composition/decomposition rule.

\subsection{Composition}\label{Sect-Comp}

\begin{figure}[h]
\begin{center}
  \includegraphics[width=0.55\textwidth]{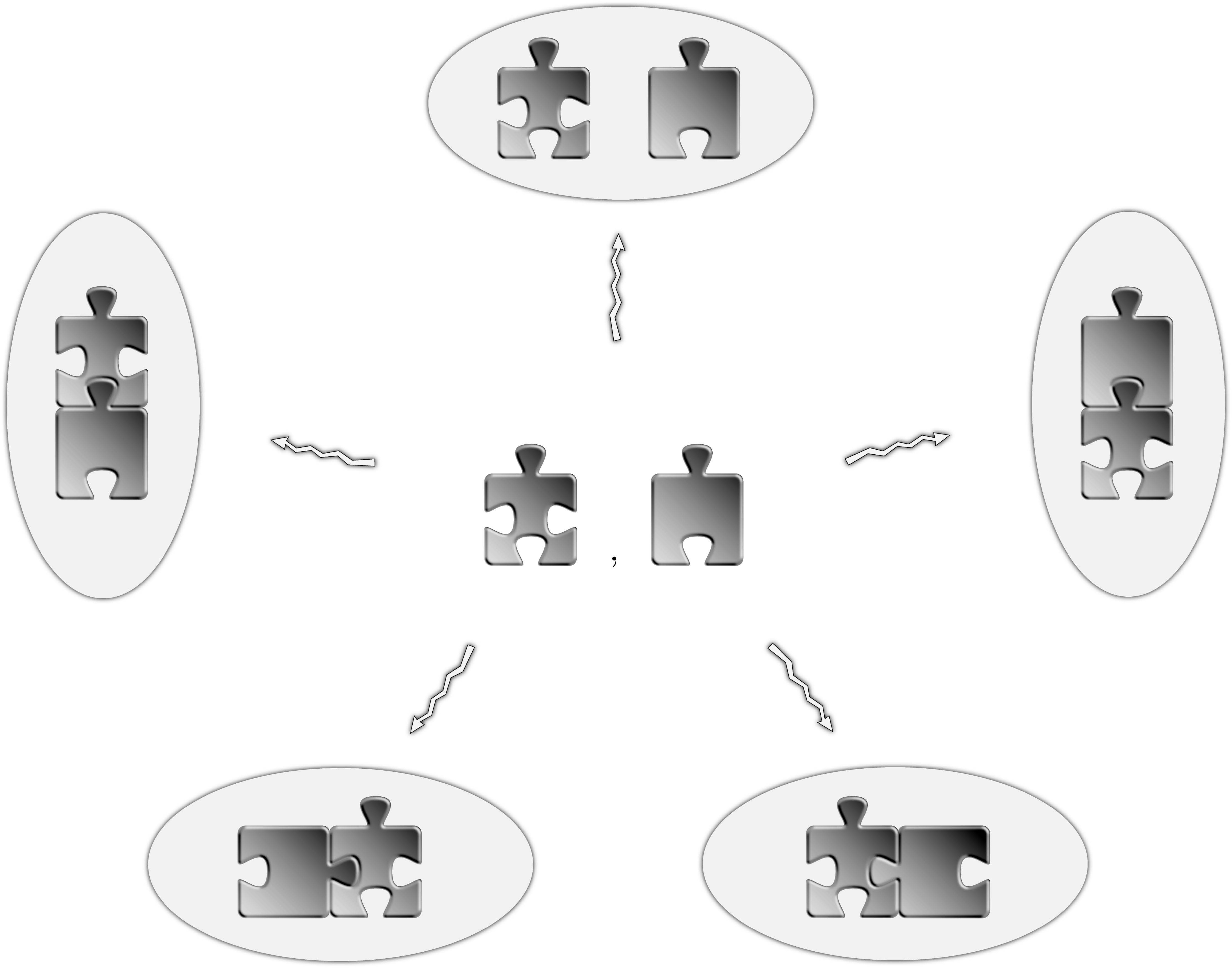}
\caption{\label{puzzle-comp} Illustration of the concept of composition. Two puzzles may compose in six possible ways, two of which give the same result (two at the bottom).}
\end{center}
\end{figure}

\emph{Composition} of objects in a combinatorial class $\mathcal{C}$ is a prescription how from two objects make another one in the same class. In general, such a rule may be indeterminate that is allow two given objects to compose in a number of ways. Furthermore, it can happen that some of these various possibilities produce the same outcome. See Fig.~\ref{puzzle-comp} for illustration. Therefore, a complete description of composition should keep a record of all the options which is conveniently attained by means of the multiset construction. Here is the formal definition:

\begin{definition}[Composition]\label{Comp}\ \vspace{0.1cm}\\
For a given combinatorial class $\mathcal{C}$ the composition rule is a mapping
\begin{eqnarray}\label{triangle}
\blacktriangleleft\ :\mathcal{C}\times\mathcal{C}&\longrightarrow&\textsc{MSet}\,(\mathcal{C})\ ,
\end{eqnarray}
assigning to each pair of objects $\varGamma_2,\varGamma_1\in\mathcal{C}$ the multiset, denoted by $\varGamma_2\blacktriangleleft\varGamma_1$, consisting of all possible compositions of $\varGamma_2$ with $\varGamma_1$, where multiple copies keep an account of the number of ways in which given outcome occurs. Sometimes, for brevity we will write  
$(\varGamma_2,\varGamma_1) \leadsto\varGamma\in \varGamma_2\blacktriangleleft\varGamma_1$.
\end{definition}
Note that this definition naturally extends to the mapping $\blacktriangleleft\ :\textsc{MSet}\,(\mathcal{C})\times\textsc{MSet}\,(\mathcal{C})\longrightarrow\textsc{MSet}\,(\mathcal{C})$ which for given two multisets $\varGamma_2,\varGamma_1\in\textsc{MSet}(\mathcal{C})$ take their elements one by one, compose and collect the results all together, \textit{i.e.} 
\begin{eqnarray}
\varGamma_2\blacktriangleleft\varGamma_1=\biguplus_{\gamma_2\in\varGamma_2,\gamma_1\in\varGamma_1}\gamma_2\blacktriangleleft\gamma_1\ .
\end{eqnarray}

At this point the concept of composition is quite general, and its further development obviously depends on the choice of the rule. One supplements this construction with additional constraints. Below we discuss some natural conditions one might expect from a reasonable composition rule.

\begin{itemize}
\item[\textbf{\texttt{(C1)}}] {\textbf{\emph{Finiteness.}}
It is sensible to assume that objects compose only in a finite number of ways, \textit{i.e.} for each $\varGamma_2,\varGamma_1\in\mathcal{C}$ we have
\begin{eqnarray}\label{Finiteness-comp}
\#\ \varGamma_2\blacktriangleleft\varGamma_1<\infty\ .
\end{eqnarray}}
\item[\textbf{\texttt{(C2)}}] {\textbf{\emph{Triple composition.}}
Composition applies to more that two objects as well. For given $\varGamma_3,\varGamma_2,\varGamma_1\in\mathcal{C}$ one can compose them successively and construct the multiset of possible compositions. There are two possible scenarios however: one can either start by composing the first two $(\varGamma_2,\varGamma_1)\leadsto\varGamma'$ and then composing the outcome with the third $(\varGamma_3,\varGamma')\leadsto\varGamma$, or change the order and begin with $(\varGamma_3,\varGamma_2)\leadsto\varGamma''$ followed by the composition with the first $(\varGamma'',\varGamma_1)\leadsto\varGamma$. It is plausible to require that both scenarios lead to the same multiset. This condition is a sort of associativity property which in a compact form reads
\begin{eqnarray}\label{Triple-comp}
\varGamma_3\blacktriangleleft(\varGamma_2\blacktriangleleft\varGamma_1)=(\varGamma_3\blacktriangleleft\varGamma_2)\blacktriangleleft\varGamma_1\ .
\end{eqnarray}
Note that it justifies dropping of the brackets in the denotation of triple composition $\varGamma_3\blacktriangleleft\varGamma_2\blacktriangleleft\varGamma_1$. Clearly, the procedure generalizes to multiple compositions and Eq.~(\ref{Triple-comp}) entails analogous condition in this case as well.}
\item[\textbf{\texttt{(C3)}}] {\textbf{\emph{Neutral object.}}
Often, in a class there exists a neutral object, denoted by \O, which composes with elements of the class only in a trivial way, \textit{i.e.} $(\text{\O},\varGamma)\leadsto\varGamma$ and $(\varGamma,\text{\O})\leadsto\varGamma$. In other words, for each $\varGamma\in\mathcal{C}$ we have
\begin{eqnarray}
\text{\O}\blacktriangleleft\varGamma=\{\varGamma\}\ \ \ \ \ \&\ \ \ \ \ \varGamma \blacktriangleleft\text{\O}=\{\varGamma\}\ .
\end{eqnarray}
Note that if \O\ exists, it is unique.}
\item[\textbf{\texttt{(C4)}}] {\textbf{\emph{Symmetry.}}
Sometimes the composition rule is such that the order in which elements are composed is irrelevant. Then for each $\varGamma_2,\varGamma_1\in\mathcal{C}$ the following commutativity condition holds
\begin{eqnarray}
\varGamma_2\blacktriangleleft\varGamma_1=\varGamma_1\blacktriangleleft\varGamma_2\ .
\end{eqnarray}}
\end{itemize}

\subsection{Decomposition}\label{Sect-Decomp}

\begin{figure}[t]
\begin{center}
  \includegraphics[width=0.75\textwidth]{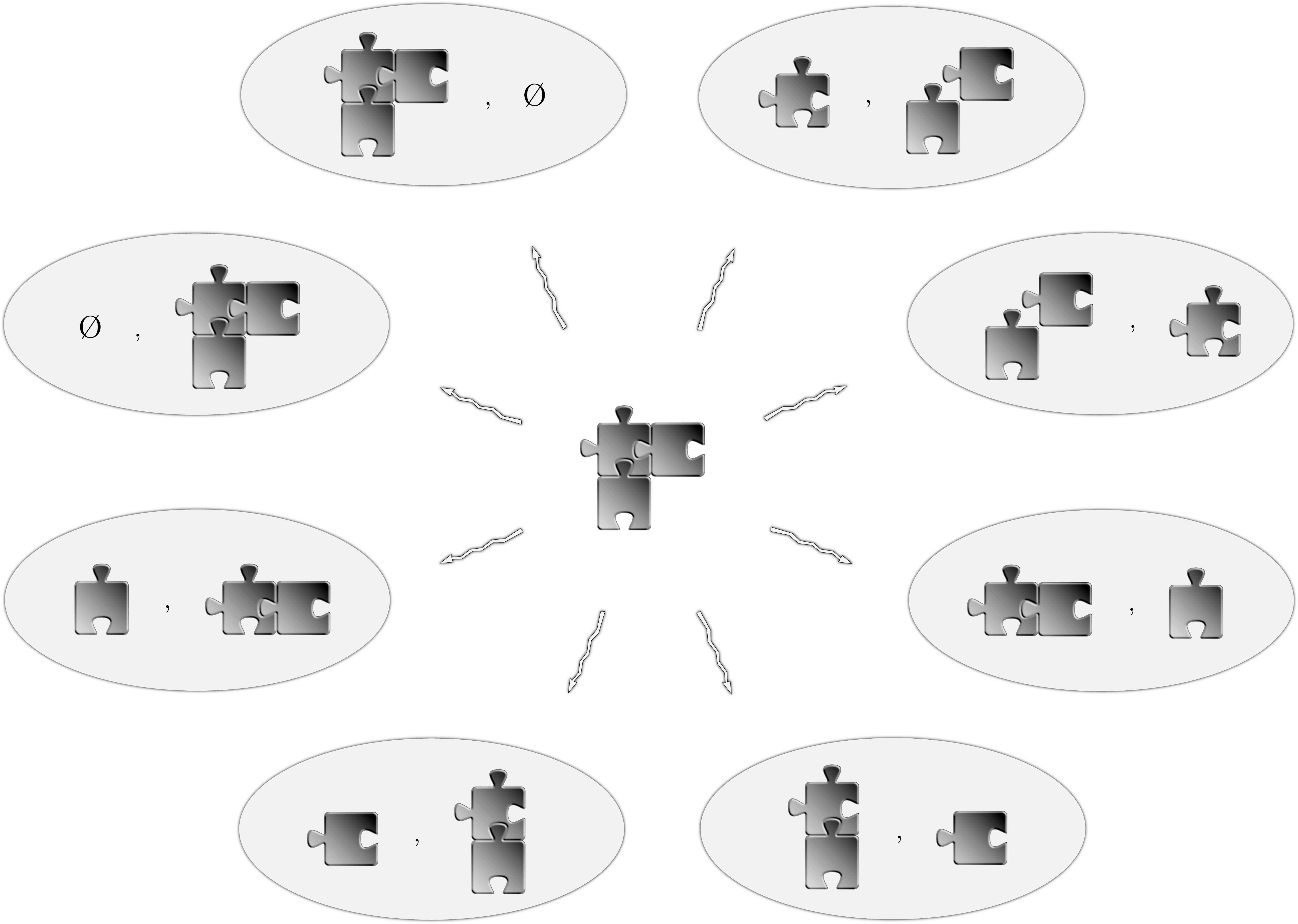}
\caption{\label{puzzle-decomp} Illustration of the concept of decomposition. A compound puzzle decomposes into ordered pairs of pieces in eight possible ways. Note that two outcomes occur twice: the two at the bottom-left corner and the two at the bottom-right are the same.}
\end{center}
\end{figure}

Suppose that a combinatorial class $\mathcal{C}$ allows for \emph{decomposition} of objects, \emph{i.e.} splitting into ordered pairs of pieces from the same class. In general, there might be various ways of splitting an object following  a given rule and, moreover, some of them may yield the same result. See Fig.~\ref{puzzle-decomp} for illustration. The whole collection of possibilities is again properly described by the notion of multiset. Hence, we have the definition

\begin{definition}[Decomposition]\label{Decomp}\ \vspace{0.1cm}\\
Decomposition rule in a combinatorial class $\mathcal{C}$ is a mapping
\begin{eqnarray}\label{angle}
\langle\,\cdot\,\rangle:\mathcal{C}&\longrightarrow&\textsc{MSet}\,(\mathcal{C}\times\mathcal{C})\ ,
\end{eqnarray}
which for each object $\varGamma \in\mathcal{C}$ defines the multiset, denoted by $\langle\varGamma\rangle$, comprised of all pairs $(\varGamma'',\varGamma')$ which are splittings of $\varGamma$, with multiple copies keeping a record of possible decompositions producing the same outcome. Concisely, we will write $\varGamma \leadsto(\varGamma'',\varGamma')\in \langle\varGamma\rangle$.
\end{definition}
Extension of the definition to the mapping $\langle\,\cdot\,\rangle:\textsc{MSet}\,(\mathcal{C})\longrightarrow\textsc{MSet}\,(\mathcal{C}\times\mathcal{C})$ is straightforwardly given by collecting all together decompositions of elements taken from $\varGamma\in\textsc{MSet}(\mathcal{C})$, \textit{i.e.}
\begin{eqnarray}
\langle\varGamma\rangle=\ \biguplus_{\gamma\in\varGamma}\ \ \langle\gamma\rangle\ .
\end{eqnarray}

Below, analogously as in Section~\ref{Sect-Comp} we consider some general conditions which one might require from a reasonable decomposition rule. We observe that most of them are in a sense dual to those discussed for the composition rule, which reflects the opposite character of both procedures. Note, however, that the decomposition rule is so far unrelated to composition -- insofar as the latter might be even not defined -- and the conditions should be treated as independent.

\begin{itemize}
\item[\textbf{\texttt{(D1)}}] {\textbf{\emph{Finiteness.}}
One may reasonably expect that objects decompose in a finite number of ways only, \textit{i.e.} for each $\varGamma\in\mathcal{C}$ we require
\begin{eqnarray}
\#\ \langle\varGamma\rangle<\infty\ .
\end{eqnarray}}
\item[\textbf{\texttt{(D2)}}] {\textbf{\emph{Triple decomposition.}}
Decomposition into pairs naturally extends to splitting  an object into three pieces $\varGamma\leadsto (\varGamma_3,\varGamma_2,\varGamma_1)$. An obvious way to carry out the multiple splitting is by applying the same procedure repeatedly, \textit{i.e.} decomposing one of the components obtained in the preceding step. Following this prescription one usually expects that the result does not depend on the choice of the component it is applied to. In other words, we require that we  end up with the same collection of triple decompositions when splitting $\varGamma\leadsto  (\varGamma'',\varGamma_1)$ and then splitting the left component $\varGamma''\leadsto  (\varGamma_3,\varGamma_2)$, 
as in the case when starting with $\varGamma\leadsto  (\varGamma_3,\varGamma')$ and then splitting the right component $\varGamma'\leadsto  (\varGamma_2,\varGamma_1)$.
This condition can be seen as a sort of co-associativity property for decomposition, and in explicit form boils down to the following equality between multisets of triple decompositions
\begin{eqnarray}\label{Triple-decomp}
\biguplus_{(\varGamma'',\varGamma')\in\langle\varGamma\rangle} \{\varGamma''\}\,\times\, \langle\varGamma'\rangle
=\biguplus_{(\varGamma'',\varGamma')\in\langle\varGamma\rangle}\langle\varGamma''\rangle\,\times\, \{\varGamma'\}\ .
\end{eqnarray}
The above procedure directly extends to splitting into multiple pieces $\varGamma\leadsto (\varGamma_n,...\varGamma_1)$ by iterated decomposition. Clearly, the condition of Eq.~(\ref{Triple-decomp}) asserts the same result no matter which way decompositions are carried out. Hence, we can consistently define the multiset consisting of multiple decompositions of an object as
\begin{eqnarray}\label{Multiple-decomp}
\langle\varGamma\rangle^{(n)}=\biguplus_{\varGamma\leadsto (\varGamma_n,...\varGamma_1)}\{(\varGamma_n,...,\varGamma_1)\}\ ,
\end{eqnarray}
with the convention $\langle\varGamma\rangle^{(1)}=\langle\varGamma\rangle$.}
\item[\textbf{\texttt{(D3)}}] {\textbf{\emph{Void object.}}
Oftentimes, a class contains a void (or empty) element \O, such that objects decompose in a trivial way. It should have the property that any object $\varGamma\neq\text{\O}$ splits into a pair containing either \O\ or $\varGamma$ in exactly two ways
\begin{eqnarray}
\varGamma\leadsto (\text{\O},\varGamma)\ \ \ \ \ \ \&\ \ \ \ \ \ \varGamma\leadsto (\varGamma,\text{\O})\ ,
\end{eqnarray}
and $\text{\O}\leadsto (\text{\O},\text{\O})$. Clearly, if \O\ exists, it is unique.}

\item[\textbf{\texttt{(D4)}}] {\textbf{\emph{Symmetry.}}
For some rules the order between components in decompositions  is immaterial, \emph{i.e.} it allows for the exchange $(\varGamma',\varGamma'')\longleftrightarrow(\varGamma'',\varGamma')$. In this case we have the following symmetry condition
\begin{eqnarray}\label{symmetry}
(\varGamma',\varGamma'')\in\langle\varGamma\rangle\ \Longleftrightarrow\ (\varGamma'',\varGamma')\in\langle\varGamma\rangle\ ,
\end{eqnarray}
and multiplicities of $(\varGamma',\varGamma'')$ and $(\varGamma'',\varGamma')$ in $\langle\varGamma\rangle$ are the same.}
\item[\textbf{\texttt{(D5)}}] {\textbf{\emph{Finiteness of multiple decompositions.}}
Recall multiple decompositions $\varGamma\leadsto(\varGamma_n,...\varGamma_1)$ considered in condition \textbf{\texttt{(D4)}} and observe that we may go with the number of components to any $n\in\mathbb{N}$.  However, if one takes into account only nontrivial decompositions, \textit{i.e.} such that do not contain void \O\ components, it is often the case that the process terminates after a finite number of steps.
In other words, for each $\varGamma\in\mathcal{C}$ there exists $N\in\mathbb{N}$ such that for all $n\geqslant N$ one has
\begin{eqnarray}\label{FinMultDecomp}
\left\{\varGamma\leadsto(\varGamma_n,...\varGamma_1): \varGamma_n,...,\varGamma_1\neq\text{\O}\right\}=\emptyset\ .
\end{eqnarray}}
\end{itemize}

\subsection{Compatibility}\label{Sect-Compatibility}

Now, let us take a combinatorial class $\mathcal{C}$ which admits both composition and decomposition of objects at the same time. We give a simple \emph{compatibility} condition for both procedures to be consistently combined together.
\begin{itemize}
\item[\textbf{\texttt{(CD1)}}] {\textbf{\emph{Composition--decomposition compatibility.}}
Suppose we are given a pair of objects $(\varGamma_2,\varGamma_1)\in\mathcal{C}\times\mathcal{C}$ which we want to decompose. We may think of two consistent decomposition schemes which involve composition as an intermediate step. We can either start by composing them together $\varGamma_2\plug{}\varGamma_1$ and then splitting all the resulting objects into pieces, or first  decompose each of them separately into $\langle\varGamma_2\rangle$ and $\langle\varGamma_1\rangle$ and then compose elements of both multisets in a component-wise manner. One may reasonably expect the same outcome no matter which procedure is applied. The formal description of compatibility comes down to the equality of multisets
\begin{eqnarray}\label{compatibility}
\langle\varGamma_2\blacktriangleleft\varGamma_1\rangle=\biguplus_{\substack{(\varGamma_2'',\varGamma_2')\in\langle\varGamma_2\rangle\\(\varGamma_1'',\varGamma_1')\in\langle\varGamma_1\rangle}}(\varGamma_2''\blacktriangleleft\varGamma_1'')\times(\varGamma_2'\blacktriangleleft\varGamma_1')\ .
\end{eqnarray}
We remark that this property implies that the void and neutral object of conditions \textbf{\texttt{(D3)}} and \textbf{\texttt{(C3)}} are the same, and hence common denotation \O.}
\end{itemize}

Oftentimes, composition/decomposition of objects come alongside with the notion of size. It is usually the case when their defining rules make use of the same characteristics which are counted by the size function. Here is a useful condition connecting these concepts: 
\begin{itemize}
\item[\textbf{\texttt{(CD2)}}] {\textbf{\emph{Compatibility with size.}} It may happen that the considered composition rule preserves size. This means that when composing two objects $(\varGamma_2,\varGamma_1)\leadsto\varGamma$ sizes of both components add up, \textit{i.e.}
\begin{eqnarray}\label{Size-comp}
|\varGamma|=|\varGamma_2|+|\varGamma_1|\ \ \ \ \ \ \text{for}\ \ \ \ \ \ \varGamma\in\varGamma_2\blacktriangleleft\varGamma_1\ .
\end{eqnarray}
This requirement boils down to the restriction on the mapping of Eq.~(\ref{triangle}) in the following way:\footnote{\label{footnote-size}Recall, that the size function $|\cdot|:\mathcal{C}\longrightarrow\naturals$ divides $\mathcal{C}$ into disjoint subclasses $\mathcal{C}_n=\{\varGamma\in\mathcal{C}:|\varGamma|=n\}$, such that $\mathcal{C}=\bigcup_{n\in\naturals}\mathcal{C}_n$.}
\begin{eqnarray}
\blacktriangleleft\ :\mathcal{C}_i\times\mathcal{C}_j&\longrightarrow&\textsc{MSet}\,(\mathcal{C}_{i+j})\ .
\end{eqnarray}
A parallel condition for decomposition implies that after splitting $\varGamma\leadsto(\varGamma'',\varGamma')$ the original size of an object distributes between the parts, \textit{i.e.}
\begin{eqnarray}\label{Size-decomp}
|\varGamma|=|\varGamma''|+|\varGamma'|\ \ \ \ \ \ \text{for}\ \ \ \ \ \ (\varGamma'',\varGamma')\in\langle\varGamma\rangle\ .
\end{eqnarray}
This translates into the constraint on the mapping of Eq.~(\ref{angle}) as:
\begin{eqnarray}
\langle\,\cdot\,\rangle:\mathcal{C}_k&\longrightarrow&\textsc{MSet}\,\Big(\biguplus_{i+j=k}\mathcal{C}_i\times\mathcal{C}_j\ \Big)\ .
\end{eqnarray}
In the following, we will assume that there is a single object of size zero, \textit{i.e.} $\mathcal{C}_0=\{\text{\O}\}$.}
\end{itemize}

\section{Construction of Algebraic Structures}\label{Sect-Algebra}

We will demonstrate the way in which combinatorial objects can be equipped with natural algebraic structures based on the composition/decomposition concept. The key role in the argument play the conditions discussed in Section~\ref{Sect-Comp-Decomp} which provide a route to systematic construction of the algebra, co-algebra, bi-algebra and Hopf algebra structures. We note that most of combinatorial algebras can be systematically treated along these lines.

\subsection{Vector space}

For a given combinatorial class $\mathcal{C}$ we will consider a vector space $\mathscr{C}$ over a field $\mathbb{K}$ which consists of (finite) linear combinations of elements in $\mathcal{C}$, \textit{i.e.} 
\begin{eqnarray}\label{VectorSpace}
\mathscr{C}=\mathbb{K}\,\mathcal{C}=\left\{\ {\sum}_i\alpha_i\ \varGamma_i:\ \alpha_i\in\mathbb{K},\ \varGamma_i\in\mathcal{C}\ \right\}\ .
\end{eqnarray}
Addition of elements and multiplication by scalars in $\mathscr{C}$ has the usual form
\begin{eqnarray}\label{Addition}
{\sum}_i\ \alpha_i\ \varGamma_i+{\sum}_i\ \beta_i\ \varGamma_i&=&{\sum}_i\ (\alpha_i+\beta_i)\ \varGamma_i\ ,\\
\alpha\ {\sum}_i\ \beta_i\ \varGamma_i&=&{\sum}_i\ (\alpha\,\beta_i)\ \varGamma_i\ .
\end{eqnarray}
Clearly, elements of $\mathcal{C}$ are independent and span the whole vector space. Hence, $\mathscr{C}$ comes endowed with the distinguished basis which, in addition, carries a combinatorial meaning. We will call $\mathcal{C}$ the \emph{combinatorial basis} of $\mathscr{C}$.

\subsection{Multiplication \& co-multiplication}\label{Sect-Mult-Co-mult}

Having defined the vector space $\mathscr{C}$ built on a combinatorial class $\mathcal{C}$ we are ready to make use of its combinatorial content. Below we provide a general scheme for constructing an algebra and co-algebra structures \cite{Bo89} based on the notions of composition and decomposition discussed in Section~\ref{Sect-Comp-Decomp}.

Suppose $\mathcal{C}$ admits \emph{composition} as defined in Section~\ref{Sect-Comp}. We will consider a bilinear mapping 
\begin{eqnarray}
*:\mathscr{C}\times\mathscr{C}\longrightarrow\mathscr{C}
\end{eqnarray}
defined on basis elements $\varGamma_2,\varGamma_1\in\mathcal{C}$ as the sum of all possible compositions of $\varGamma_2$ with $\varGamma_1$, \textit{i.e.}
\begin{eqnarray}\label{multiplication}
\varGamma_2*\varGamma_1=\sum_{\varGamma\in\varGamma_2\blacktriangleleft\varGamma_1}\varGamma\ .
\end{eqnarray}
Note, that although all coefficients in the defining Eq.~(\ref{multiplication}) are equal to one, some of the terms in the sum may appear several times; this is because $\varGamma_2\blacktriangleleft\varGamma_1$ is a multiset.  One rightly anticipates that multiplicities of elements will play the role of \emph{structure constants} of the algebra.
Such defined mapping is a natural candidate for \emph{multiplication} and we have the following statement

\begin{proposition}[Algebra] \label{Prop-Algebra} \ \vspace{0.1cm}\\
The vector space $\mathscr{C}$ with the multiplication defined in Eq.~(\ref{multiplication}) forms an associative algebra with unit $(\mathscr{C},+,*,\emph{\text{\O}})$ if conditions \emph{\textbf{\texttt{(C1)}}} -- \emph{\textbf{\texttt{(C3)}}} hold. Under condition \emph{\textbf{\texttt{(C4)}}} it is commutative.
\end{proposition}

\begin{proof}
Condition \textbf{\texttt{(C1)}} guarantees that the sum in Eq.~(\ref{multiplication}) is finite -- hence it is well defined. Conditions \textbf{\texttt{(C3)}} and \textbf{\texttt{(C4)}} directly translate into the existence of the unit element  \O\ and commutativity respectively.
Associativity is the consequence of bilinearity of multiplication and condition \textbf{\texttt{(C2)}} which asserts equality of multisets resulting from two scenarios of triple composition $(\varGamma_3,\varGamma_2,\varGamma_1)\leadsto \varGamma$; it is straightforward to check for basis elements that
\begin{eqnarray}
\varGamma_3*(\varGamma_2*\varGamma_1)=\sum_{\varGamma\in\varGamma_3\blacktriangleleft\varGamma_2\blacktriangleleft\varGamma_1}\varGamma\ \ =(\varGamma_3*\varGamma_2)*\varGamma_1\ .
\end{eqnarray}

\end{proof}

Now, we will consider $\mathcal{C}$ equipped with the notion of \emph{decomposition} as described in Section~\ref{Sect-Decomp}. Let us take a linear mapping 
\begin{eqnarray}
\Delta:\mathscr{C}&\longrightarrow&\mathscr{C}\otimes\mathscr{C}
\end{eqnarray}
defined on basis elements $\varGamma\in\mathcal{C}$ as the sum of all splittings into pairs, which in explicit form reads
\begin{eqnarray}\label{delta}
\Delta\,(\varGamma)=\sum_{(\varGamma'',\varGamma')\in\langle\varGamma\rangle}\varGamma''\otimes\varGamma'\ .
\end{eqnarray}
Repetition of terms in Eq.~(\ref{delta}) leads after simplification to coefficients which are multiplicities of elements in the multiset $\langle\varGamma\rangle$. These numbers are sometimes called \emph{section coefficients}, see \cite{JoniRota}.
We will also need a linear mapping 
\begin{eqnarray}
\varepsilon:\mathscr{C}&\longrightarrow&\mathbb{K}
\end{eqnarray}
which extracts the expansion coefficient standing at the void \O. It is defined on basis elements $\varGamma\in\mathcal{C}$ in a canonical way
\begin{eqnarray}\label{epsilon}
\varepsilon(\varGamma)=\left\{
\begin{array}{l}
1 \text{\ \ \ \ if\ \ \ }\varGamma=\text{\O}\ ,\\
0  \text{\ \ \ \ otherwise\ . }
\end{array}\right.
\end{eqnarray}
These mappings play the role of \emph{co-multiplication} and \emph{co-unit} in the construction of a co-algebra as explained in the following proposition

\begin{proposition}[Co-algebra] \label{Prop-Co-algebra} \ \vspace{0.1cm}\\
If conditions \emph{\textbf{\texttt{(D1)}}} -- \emph{\textbf{\texttt{(D3)}}} are satisfied the mappings $\Delta$ and $\varepsilon$ defined in Eqs.~(\ref{delta}) and (\ref{epsilon}) respectively are the co-multiplication and co-unit which make the vector space $\mathscr{C}$ into a co-algebra $(\mathscr{C},+,\Delta,\varepsilon)$. It is co-commutative if condition \emph{\textbf{\texttt{(D4)}}} holds.
\end{proposition}

\begin{proof}
The sum in Eq.~(\ref{delta}) is well defined as long as the number of compositions is finite, \textit{i.e.} condition \textbf{\texttt{(D1)}} is satisfied. 
From equivalence of triple splittings $\varGamma\leadsto(\varGamma_3,\varGamma_2,\varGamma_1)$ obtained in two possible ways considered in condition \textbf{\texttt{(D2)}}, one readily verifies for a basis element that
\begin{eqnarray}
(\Delta\otimes Id)\circ\Delta\, (\varGamma)=\sum_{\varGamma\leadsto(\varGamma_3,\varGamma_2,\varGamma_1)}\varGamma_3\otimes\varGamma_2\otimes\varGamma_1\ =(Id\otimes\Delta)\circ\Delta\, (\varGamma)\ ,
\end{eqnarray}
which by linearity extends on all $\mathscr{C}$ proving co-associativity of co-multiplication defined in Eq.~(\ref{delta}). \\
The co-unit $\varepsilon: \mathscr{C}\longrightarrow \mathbb{K}$ by definition should satisfy the equalities
\begin{eqnarray}\label{co-unit-def}
(\varepsilon\otimes Id)\circ\Delta=Id=(Id\otimes\varepsilon)\circ\Delta \ ,
\end{eqnarray}
where the identification $\mathbb{K}\otimes\mathscr{C}=\mathscr{C}\otimes\mathbb{K}=\mathscr{C}$ is implied.
We will check the first one for a basis element $\varGamma$ by direct calculation
\begin{eqnarray}
(\varepsilon\otimes Id)\circ\Delta\,(\varGamma)
=\sum_{(\varGamma_1,\varGamma_2)\in\langle\varGamma\rangle}\varepsilon(\varGamma_1)\otimes\varGamma_2=1\otimes\varGamma=\varGamma=Id\,(\varGamma)\ .
\end{eqnarray}
Note that we have applied condition \textbf{\texttt{(D3)}} by taking all terms in the sum equal to zero except the unique decomposition $(\text{\O},\varGamma)$ picked up by $\varepsilon$ in accordance with the definition of Eq.~(\ref{epsilon}). The identification $1\otimes\varGamma=\varGamma$ completes the proof of the first equality in Eq.~(\ref{co-unit-def}); verification of the second one is analogous.
\\
Check of co-commutativity of the co-product under condition \textbf{\texttt{(D4)}} is immediate.

\end{proof}

\subsection{Bi-algebra and Hopf algebra structure}\label{Sect-Bi-algebra}
We have seen in Propositions~\ref{Prop-Algebra} and \ref{Prop-Co-algebra} how the notions of composition and decomposition lead to algebra and co-algebra structure respectively. Both schemes can be combined together so to furnish $\mathscr{C}$ with a bi-algebra structure.

\begin{theorem}[Bi-algebra] \label{bi-algebra} \ \vspace{0.1cm}\\
If condition \emph{\textbf{\texttt{(CD1)}}} is satisfied, than the algebra and co-algebra structures in $\mathscr{C}$ are compatible and combine into a bi-algebra $(\mathscr{C},+,*,\emph{\text{\O}},\Delta,\varepsilon)$.
\end{theorem}

\begin{proof}
The structure of a bi-algebra requires that the co-multiplication $\Delta:\mathscr{C}\otimes\mathscr{C}\longrightarrow\mathscr{C}$ and the co-unit $\varepsilon:\mathscr{C}\longrightarrow\mathbb{K}$ of the co-algebra preserve multiplication in $\mathscr{C}$. Thus, we need to verify for basis elements $\varGamma_1$ and $\varGamma_2$ that
\begin{eqnarray}\label{delta-morph}
\Delta\,(\varGamma_2*\varGamma_1)&=&\Delta\,(\varGamma_2)*\Delta\,(\varGamma_1)\ ,
\end{eqnarray}
with component-wise multiplication in the tensor product $\mathcal{G}\otimes\mathcal{G}$ on the right-hand-side, and that
\begin{eqnarray}\label{epsilon-morph}
\varepsilon\,(\varGamma_2*\varGamma_1)&=&\varepsilon\,(\varGamma_2)\ \varepsilon\,(\varGamma_1)\ ,
\end{eqnarray}
with terms on the right-hand-side multiplied in the field $\mathbb{K}$.
\\
We check Eq.~(\ref{delta-morph}) directly by expanding both  sides using definitions of Eqs.~(\ref{multiplication}) and (\ref{delta}). Accordingly, the left-hand-side takes the form
\begin{eqnarray}\label{DM1}
\Delta\,(\varGamma_2*\varGamma_1)=\sum_{\varGamma\in\varGamma_2\plug{}\varGamma_1}\Delta\,(\varGamma)
=\sum_{\substack{\varGamma\in\varGamma_2\plug{}\varGamma_1\\(\varGamma'',\varGamma')\in\langle\varGamma\rangle}}\varGamma''\otimes\varGamma'=\sum_{(\varGamma'',\varGamma')\in\langle\varGamma_2\plug{}\varGamma_1\rangle}\varGamma''\otimes\varGamma'\ ,
\end{eqnarray}
while the right-hand-side reads
\begin{eqnarray}\label{DM2}
\Delta\,(\varGamma_2)*\Delta\,(\varGamma_1)=\sum_{\substack{(\varGamma_2'',\varGamma_2')\in\langle\varGamma_2\rangle\\(\varGamma_1'',\varGamma_1')\in\langle\varGamma_1\rangle}}\underbrace{(\varGamma_2''\otimes\varGamma_2')*(\varGamma_1''\otimes\varGamma_1')}_{(\varGamma_2''*\varGamma_1'')\otimes(\varGamma_2'*\varGamma_1')}
=\sum_{\substack{(\varGamma_2'',\varGamma_2')\in\langle\varGamma_2\rangle\\(\varGamma_1'',\varGamma_1')\in\langle\varGamma_1\rangle}}\ \ \sum_{\substack{\varGamma''\in\varGamma_2''\plug{}\varGamma_1''\\\varGamma'\in\varGamma_2'\plug{}\varGamma_1'}}\varGamma''\otimes\varGamma'\ .
\end{eqnarray}
A closer look at condition \textbf{\texttt{(CD1)}} and Eq.~(\ref{compatibility}) shows a one-to-one correspondence between terms in the sums on the right-hand-sides of Eqs.~(\ref{DM1}) and (\ref{DM2}), which proves Eq.~(\ref{delta-morph}).
\\
Verification of Eq.~(\ref{epsilon-morph}) rests upon a simple observation, steaming from \textbf{\texttt{(C3)}}, \textbf{\texttt{(D3)}} and \textbf{\texttt{(CD1)}}, that composition of objects $\varGamma_2\plug{}\varGamma_1$ yields the neutral element \O\ only if both of them are void. Then, both sides of Eq.~(\ref{epsilon-morph}) are equal to $1$ for $\varGamma_1=\varGamma_2=\text{\O}$ and $0$ otherwise, which ends the proof.

\end{proof}

Finally, let us take a linear mapping
\begin{eqnarray}\label{antipode}
S:\mathscr{C}&\longrightarrow&\mathscr{C}\ ,
\end{eqnarray}
defied as an alternating sum of multiple products over possible nontrivial decompositions of an object, \textit{i.e.}
\begin{eqnarray}\label{antipode}
S(\varGamma)=\sum_{\substack{\varGamma\leadsto(\varGamma_n,...,\varGamma_1)\\\varGamma_n,...,\varGamma_1\neq\,\text{\O}}}(-1)^n\ \varGamma_n*...*\varGamma_1
\end{eqnarray}
for $\varGamma\neq\text{\O}$ and $S(\text{\O})=\text{\O}$. This mapping provides an \emph{antipode} completing the construction of a Hopf algebra structure on $\mathscr{C}$ \cite{SweedlerBook,Ab80}.

\begin{proposition}[Hopf Algebra] \label{HopfAlgebra} \ \vspace{0.1cm}\\
If furthermore condition \emph{\textbf{\texttt{(D5)}}} holds, $S$ defined in Eq.~(\ref{antipode}) is the antipode which makes the bi-algebra  $\mathscr{C}$ of Theorem~\ref{bi-algebra} into a Hopf algebra $(\mathscr{C},+,*,\emph{\text{\O}},\Delta,\varepsilon,S)$.
\end{proposition}

\begin{proof}

A Hopf algebra consists of a bi-algebra $(\mathscr{C},+,*,\text{\O},\Delta,\varepsilon)$ equipped with an antipode $S:\mathscr{C}\longrightarrow\mathscr{C}$. The latter is an endomorphism by definition satisfying the property
\begin{eqnarray}\label{antipode-def}
\mu\circ(Id\otimes S)\circ\Delta=\epsilon=\mu\circ(S\otimes Id)\circ\Delta\ ,
\end{eqnarray}
where for better clarity multiplication was denoted by $\mu(\varGamma_2\otimes \varGamma_1)=\varGamma_2*\varGamma_1$. The mapping $\epsilon:\mathscr{C}\longrightarrow\mathscr{C}$, defined as $\epsilon=\text{\O}\,\varepsilon$, is the projection on the subspace spanned by \O, \emph{i.e.}
\begin{eqnarray}\label{epsilon-proj}
\epsilon(\varGamma)=\left\{
\begin{array}{l}
\varGamma \text{\ \ \ \ if\ \ \ }\varGamma=\text{\O}\ ,\\
0  \text{\ \ \ \ otherwise\ . }
\end{array}\right.
\end{eqnarray}
We will prove that $S$ given in Eq.~(\ref{antipode}) satisfies the condition of Eq.~(\ref{antipode-def}). We start by considering an auxiliary linear mapping $\Phi:\mathsf{End}(\mathscr{C})\longrightarrow \mathsf{End}(\mathscr{C})$ defined as
\begin{eqnarray}\label{Phi}
\Phi(f)=\mu\circ(Id\otimes f)\circ\Delta\ ,\ \ \ \ \ \ \ \text{for} \ \ \ f\in \mathsf{End}(\mathscr{C})\ .
\end{eqnarray}
Observe that under the assumption that $\Phi$ is invertible the first equality in Eq.~(\ref{antipode-def}) can be rephrased into the condition
\begin{eqnarray}\label{antipode-inverse-Phi}
S=\Phi^{-1}(\epsilon)\ .
\end{eqnarray}
Now, our objective is to show that $\Phi$ is invertible and calculate its inverse explicitly.
By extracting identity we get $\Phi=Id+\Phi^+$ and observe that such defined $\Phi^+$ can be written in the form
\begin{eqnarray}
\Phi^+(f)=\mu\circ(\bar{\epsilon}\otimes f)\circ\Delta\ ,\ \ \ \ \ \ \ \text{for} \ \ \ f\in \mathsf{End}(\mathscr{C})\ ,
\end{eqnarray}
where $\bar{\epsilon}=Id-\epsilon$ is the complement of $\epsilon$ projecting on the subspace spanned by $\varGamma\neq\text{\O}$, \emph{i.e.}
\begin{eqnarray}\label{epsilon-proj-complement}
\bar{\epsilon}(\varGamma)=\left\{
\begin{array}{l}
0 \text{\ \ \ \ if\ \ \ }\varGamma=\text{\O}\ ,\\
\varGamma  \text{\ \ \ \ otherwise\ . }
\end{array}\right.
\end{eqnarray}
We claim that he mapping $\Phi$ is invertible with the inverse given by\footnote{For a linear mapping $L=Id+L^+:V\longrightarrow V$ its inverse can be constructed as $L^{-1}=\sum_{n=0}^\infty (-L^+)^n$ provided the sum is well defined. Indeed, one readily checks that $L\circ L^{-1}=(Id+L^+)\circ\sum_{n=0}^\infty (-L^+)^n=\sum_{n=0}^\infty (-L^+)^n+\sum_{n=0}^\infty (-L^+)^{n+1}=Id$, and similarly $L^{-1}\circ L=Id$.} 
\begin{eqnarray}\label{Phi-1}
\Phi^{-1}=\sum_{n=0}^\infty \ (-\Phi^+)^n\ .
\end{eqnarray}
In order to check that the above sum is well defined one analyzes the sum term by term. It is not difficult to calculate $n$-th iteration of $\Phi^+$ explicitly
\begin{eqnarray}\label{Phi+f}
\left(\Phi^+\right)^n(f)(\varGamma)=
\sum_{\substack{
\varGamma\leadsto(\varGamma_k,...,\varGamma_1,\varGamma_0)\\ 
\varGamma_k,...,\varGamma_1\neq\text{\O}}}
\varGamma_k*...*\varGamma_1*f(\varGamma_0)\ .
\end{eqnarray}
We note that in the above formula products of multiple decompositions arise from repeated use of the property of Eq.~(\ref{delta-morph}); the exclusion of empty components in the decompositions (except the single one on the right hand side) comes from the definition of $\bar{\epsilon}$ in Eq.~(\ref{epsilon-proj-complement}). The latter constraint together with condition \textbf{\texttt{(D5)}} asserts that the number of non-vanishing terms in Eq.~(\ref{Phi-1}) is always finite proving that $\Phi^{-1}$ is well defined.
Finally, using Eqs.~(\ref{Phi-1}) and (\ref{Phi+f}) one explicitly calculates $S$ from Eq.~(\ref{antipode-inverse-Phi}) obtaining the formula of Eq.~(\ref{antipode}).
\\
In conclusion, by construction the linear mapping $S$ of Eq.~(\ref{antipode}) satisfies the first equality in Eq.~(\ref{antipode-def}); the second equality can be checked analogously. Therefore we have proved $S$ to be an antipode thus making $\mathscr{C}$ into a Hopf algebra.

\end{proof}

We remark that by a general theory of Hopf algebras, see \cite{SweedlerBook,Ab80}, the property of Eq.~(\ref{antipode-def}) implies that $S$ is an anti-morphism and it is unique. Moreover, if  $\mathscr{C}$ is \emph{commutative} or \emph{co-commutative} $S$ is an involution, \textit{i.e.} $S\circ S=Id$.
We should also observe that the definition of the antipode given in Eq.~(\ref{antipode}) admits construction by iteration
\begin{eqnarray}
S(\varGamma)=-\sum_{\substack{(\varGamma'',\varGamma')\in\langle\varGamma\rangle\\\varGamma'\neq\ \text{\O}}}S(\varGamma'')*\varGamma'\ ,
\end{eqnarray}
and $S(\text{\O})=\text{\O}$.

Finally, whenever composition/decomposition is compatible with the notion of size in class $\mathcal{C}$ we have a \emph{grading} in the algebra $\mathscr{C}$ as explained in the following proposition:
\begin{proposition}[Grading] \label{Grading} \ \vspace{0.1cm}\\
Suppose we have a bi-algebra structure $(\mathscr{C},+,*,\emph{\text{\O}},\Delta,\varepsilon)$ constructed as in Theorem~\ref{bi-algebra}. If condition \emph{\textbf{\texttt{(CD2)}}} holds, then $\mathscr{C}$ is a graded Hopf algebra with grading given by size in $\mathcal{C}$, \textit{i.e.}
\begin{eqnarray}
\mathscr{C}=\bigoplus_{n\in\naturals}\mathscr{C}_n\ ,\ \ \ \ \ \ \ \ \ \ \ \mathscr{C}_n=\mathsf{Span}\,(\mathcal{C}_n)\ ,
\end{eqnarray}
where $\mathcal{C}_n=\{\varGamma\in\mathcal{C}:|\varGamma|=n\}$, and
\begin{eqnarray}
*:\mathscr{C}_i\times\mathscr{C}_j\longrightarrow\mathscr{C}_{i+j}\ ,\ \ \ \ \ \ \ \ \ \Delta:\mathscr{C}_k\longrightarrow \bigoplus_{i+j=k}\mathscr{C}_i\otimes\mathscr{C}_j \ .
\end{eqnarray}
\end{proposition}
\begin{proof}
Note, that condition \textbf{\texttt{(CD2)}} implies \textbf{\texttt{(D5)}}, and hence $\mathscr{C}$ is a Hopf algebra by Theorem~\ref{HopfAlgebra}. Furthermore, condition \textbf{\texttt{(CD2)}} asserts a proper action of $*$ and $\Delta$ on subspaces $\mathscr{C}_n$ built of objects of the same size.

\end{proof}

\subsection{A special case: Monoid}\label{Sect-Monoid}

Let us consider a simplified situation by taking a determinate composition law of the form
\begin{eqnarray}
\blacktriangleleft\ :\mathcal{C}\times\mathcal{C}&\longrightarrow&\mathcal{C}\ ,
\end{eqnarray}
which means that objects compose in a unique way. In other words, for each $\varGamma_2,\varGamma_1\in\mathcal{C}$ the multiset $\varGamma_2\blacktriangleleft\varGamma_1$ of Definition~\ref{Comp} is always a singleton. Observe that conditions \textbf{\texttt{(C2)}} and \textbf{\texttt{(C3)}} are equivalent to the requirement that $(\mathcal{C},\blacktriangleleft)$ is a \emph{monoid}. We note that the case of a commutative monoid was thoroughly investigated by S.A. Joni and G.-C. Rota in \cite{JoniRota} and further developed by A. Joyal \cite{Jo81}.

In this context it is convenient to consider collections $\mathfrak{C}\subset\mathcal{C}$ such that each element of $\mathcal{C}$ can be constructed as a composition of a finite number of elements from $\mathfrak{C}$, \textit{i.e.}
\begin{eqnarray}
\mathcal{C}=\left\{\ \varGamma_n\blacktriangleleft...\blacktriangleleft\varGamma_1\ :\varGamma_n,...\,,\varGamma_1\in\mathfrak{C}\ \right\}\ .
\end{eqnarray}
We call $\mathfrak{C}$ a \emph{generating class} if it is the smallest (in the sense of inclusion) subclass of $\mathcal{C}$ with this property. It has the advantage that when establishing a decomposition rule satisfying \textbf{\texttt{(CD1)}} one can specify it on the generating class $\mathfrak{C}$ in arbitrary way
\begin{eqnarray}
\langle\,\cdot\,\rangle:\mathfrak{C}&\longrightarrow&\textsc{MSet}\,(\mathcal{C}\times\mathcal{C})\ ,
\end{eqnarray}
and then consistently extend it using Eq.~(\ref{compatibility}) to the whole class $\mathcal{C}$ by defining
\begin{eqnarray}\label{extension-decomp}
\langle\varGamma_n\blacktriangleleft...\blacktriangleleft\varGamma_1\rangle\ =\biguplus_{\substack{(\varGamma_n'',\varGamma_n')\in\langle\varGamma_n\rangle\\_{\ _{\ }}.\,.\,.^{\ }\\(\varGamma_1'',\varGamma_1')\in\langle\varGamma_1\rangle}}\left\{(\varGamma_n''\blacktriangleleft...\blacktriangleleft\varGamma_1''\,,\,\varGamma_n'\blacktriangleleft...\blacktriangleleft\varGamma_1')\right\}\ .
\end{eqnarray}
We note that from a practical point of view this way of introducing the decomposition rule is very convenient as it restricts the the number of objects to be scrutinized to a smaller class $\mathfrak{C}$ and automatically guarantees compatibility of composition and decomposition rules, \textit{i.e.} \textbf{\texttt{(CD1)}} is satisfied by construction. Moreover, inspection of other properties is usually simpler in this context as well.
For example, if composition preserves size Eq.~(\ref{Size-comp}), then it is enough to check Eq.~(\ref{Size-decomp}) on $\mathfrak{C}$ and condition \textbf{\texttt{(CD2)}} automatically holds on the whole $\mathcal{C}$.

There is a canonical way in which decomposition can be introduced in this setting. Namely, one can define it on the generating elements $\varGamma\in\mathfrak{C}$ in a \emph{primitive} way, \textit{i.e.}
\begin{eqnarray}\label{primitive}
\langle\varGamma\rangle=\{(\text{\O},\varGamma), (\varGamma,\text{\O})\}\ .
\end{eqnarray}
Observe that such defined decomposition rule upon extension via Eq.~(\ref{extension-decomp}) satisfies all the conditions \textbf{\texttt{(D1)}} -- \textbf{\texttt{(D5)}}, which clears the way to the construction of a Hopf algebra. We note that objects having the property of Eq.~(\ref{primitive}) are usually called \emph{primitive elements}, for which we have
\begin{eqnarray}
\Delta(\varGamma)&=&\text{\O}\otimes\varGamma+\varGamma\otimes\text{\O}\ ,
\\
S(\varGamma)&=&-\varGamma\ .
\end{eqnarray}

\section{Examples}\label{Sect-Examples}

Here, we will illustrate how the general framework developed in Sections~\ref{Sect-Comp-Decomp} and \ref{Sect-Algebra} works in practice. We give a few examples of combinatorial objects which via composition/decomposition scheme lead to Hopf algebra structures. For other examples see \textit{e.g.} \cite{JoniRota,CartierHopfPrimer,GrLa89,BlDuHoPeSo10}.


\subsection{Words}\label{Sect-Alg-words}

Let $\mathfrak{A}=\{l_1,l_2,...,l_n\}$ be a finite set of \emph{letters} -- an \emph{alphabet}. We will consider a combinatorial class $\mathcal{A}$ consisting of (finite) words built of the alphabet $\mathfrak{A}$, \textit{i.e.} $\mathcal{A}=\mathfrak{A}^*=\{\,\text{\O},l_1,l_1l_1, l_1l_2,...\, ,l_{i_1}...\,l_{i_k},...\,\}$ where \O\ is an empty word. Size of a word will be defined as its lengths (number of letters): $|\,l_{i_1}...\,l_{i_k}|=k$ and $|\text{\O}|=0$. 
Algebraic structure in $\mathcal{A}$ can be introduced in a few ways as explained below \cite{LothaireBook,ReutenauerBook}.

\subsubsection{Free algebra of words}\label{Sect-Free-alg}
The simplest composition rule for words is given by concatenation, \textit{i.e.}\footnote{We adopt the convention that a sequence of letters indexed by the empty set is the empty word \O.}
\begin{eqnarray}\label{Free1}
l_{i_1}...\,l_{i_m}\!\blacktriangleleft\,l_{j_1}...\,l_{j_n}=l_{i_1}...\,l_{i_m}l_{j_1}...\,l_{j_n}\ .
\end{eqnarray}
Observe that $(\mathcal{A},\blacktriangleleft)$ is a monoid and $\mathfrak{A}$ is a generating class. We define decomposition of generators (letters) in the primitive way, \textit{i.e.}
\begin{eqnarray}
\langle l_i\rangle=\{(\text{\O},l_i),(l_i,\text{\O})\}\ ,
\end{eqnarray}
and extend it to the whole class $\mathcal{A}$ using Eq.~(\ref{extension-decomp}). One checks that each decomposition of a word comes down to the choice of a subword which gives the first component of a splitting (the reminder constitutes the second one), \textit{i.e.}
\begin{eqnarray}\label{Free2}
\langle l_{i_1}...\,l_{i_k}\rangle=\biguplus_{\substack{j_1<...<j_m\\j_{m+1}<...<j_{k}}}\{(l_{i_{j_1}}...\,l_{i_{j_m}},l_{i_{j_{m+1}}}...\,l_{i_{j_k}})\}
\end{eqnarray}
Note such defined composition/decomposition rule is compatible with size and hence condition \textbf{\texttt{(CD2)}} holds. Application of the scheme discussed in previous sections provide us with the mappings
\begin{eqnarray}
l_{i_1}...\,l_{i_m}*\ l_{j_1}...\,l_{j_n}=l_{i_1}...\,l_{i_m}l_{j_1}...\,l_{j_n}\ ,
\end{eqnarray}
\begin{eqnarray}\label{Free3}
\Delta( l_{i_1}...\,l_{i_k})=\sum_{\substack{j_1<...<j_m\\j_{m+1}<...<j_{k}}}l_{i_{j_1}}...\,l_{i_{j_m}}\otimes\ l_{i_{j_{m+1}}}...\,l_{i_{j_k}}\ ,
\end{eqnarray}
\begin{eqnarray}
\varepsilon(l_{i_1}...\,l_{i_k})=0\ ,\ \ \ \ \ \ \ \ \varepsilon(\text{\O})=1\ ,
\end{eqnarray}
\begin{eqnarray}\label{Free4}
S(l_{i_1}...\,l_{i_k})=(-1)^k\ l_{i_k}...\,l_{i_1}\ ,
\end{eqnarray}
which make $\mathscr{A}$ into a graded co-commutative Hopf algebra. It is called a \emph{free algebra}. Note that if the alphabet consists of more than one letter then multiplication is non-commutative.

In conclusion, we observe that if the alphabet consists of one letter only $\mathfrak{A}=\{x\}$, then the construction starts from the the class of words $\mathcal{P}=\{\text{\O},x,xx,xxx,...\}$ and leads to the \emph{algebra of polynomials} in one variable $\mathscr{P}=\mathbb{K}[x]=\left\{\,\sum_{i=0}^n \alpha_i\,x^i:\alpha_i\in\mathbb{K}\,\right\}$. In this case, we have
\begin{eqnarray}
x^i*x^j=x^{i+j}\ ,\ \ \ \ \Delta (x^n)=\sum_{i=0}^n\binom{n}{i}\  x^i\otimes x^{n-i}\ ,\ \ \ \ \varepsilon(x^n)=\delta_{n,0}\ ,\ \ \ \ S(x^n)=(-1)^n\ x^n\ .
\end{eqnarray}

\subsubsection{Symmetric algebra}

Now, let an alphabet $\mathfrak{A}=\{ l_1,l_2,...,l_n\}$ be endowed with a linear order $l_1<l_2<...<l_n$. We will consider words arranged in a non-decreasing order and define the pertaining class $\mathcal{S}$ in the form
\begin{eqnarray}
\mathcal{S}=\left\{l_{i_1}...\ l_{i_n}:i_1\leqslant...\leqslant i_n\in\naturals\right\}\ .
\end{eqnarray}
In this case, simple concatenation of words is not a legitimate composition rule and one has to amend it by additional reordering of letters
\begin{eqnarray}\label{Sym1}
l_{i_1}...\ l_{i_m}\blacktriangleleft\ l_{i_{m+1}}...\ l_{i_{m+n}}=\ l_{i_{\sigma(1)}}...\ l_{i_{\sigma(m+n)}}\ ,
\end{eqnarray}
where $\sigma$ is a unique permutation of $\{1,2,...,m+n\}$ such that $i_{\sigma(1)}\leqslant...\leqslant i_{\sigma(m+n)}$. Clearly, $(\mathcal{S},\blacktriangleleft)$ is a monoid generated by $\mathfrak{A}$. The simplest choice of primitive decomposition for the generators $\langle l_i\rangle=\{(\text{\O},l_i),(l_i,\text{\O})\}$ extends to the whole class as follows
\begin{eqnarray}\label{Sym2}
\langle l_{i_1}...\,l_{i_k}\rangle=\biguplus_{\substack{j_1<...<j_m\\j_{m+1}<...<j_{k}}}\{(l_{i_{j_1}}...\,l_{i_{j_m}},l_{i_{j_{m+1}}}...\,l_{i_{j_k}})\}
\end{eqnarray}
Observe apparent similarity of Eqs.~(\ref{Free1}) and (\ref{Free2}) to Eqs.~(\ref{Sym1}) and (\ref{Sym2}), with the only difference that words in the latter two are ordered. Construction of the \emph{symmetric algebra} $\mathscr{S}$ follows the proposed scheme, and the mappings
\begin{eqnarray}\label{SymMult}
l_n^{i_n}...\ l_1^{i_1}*\, l_n^{j_n}...\ l_1^{j_1}=\ l_n^{i_n+j_n}...\ l_1^{i_1+j_1}\ ,
\end{eqnarray}
\begin{eqnarray}\label{Free3}
\Delta( l_{i_1}...\,l_{i_k})=\sum_{\substack{j_1<...<j_m\\j_{m+1}<...<j_{k}}}l_{i_{j_1}}...\,l_{i_{j_m}}\otimes\ l_{i_{j_{m+1}}}...\,l_{i_{j_k}}\ ,
\end{eqnarray}
\begin{eqnarray}
\varepsilon(l_{i_1}...\,l_{i_k})=0\ ,\ \ \ \ \ \ \ \ \varepsilon(\text{\O})=1\ ,
\end{eqnarray}
\begin{eqnarray}\label{Free4}
S(l_{i_1}...\,l_{i_k})=(-1)^k\ l_{i_1}...\,l_{i_k}\ ,
\end{eqnarray}
define a graded Hopf algebra structure which is both commutative and co-commutative. Note that in Eq.~(\ref{SymMult}) the repeating letters were grouped together and denoted as powers. 
As a byproduct of this notation one immediately observes that the symmetric algebra $\mathscr{S}$ is isomorphic to the \emph{algebra of polynomials} in many commuting variables $\mathbb{K}[x_1,...,x_n]$.

\subsubsection{Shuffle algebra}

We will consider the class of words $\mathcal{A}$ and define composition as any \emph{shuffle} which mixes letters of the words preserving their relative order. For example, for two words ''shuffle'' and ''mix'': ''shmiufxfle'' and ''mixshuffle'' are allowed compositions whilst ''shufflemxi'' is not. Note that there are always several possible shuffles for given two (nonempty) words, and hence the use of multiset construction in definition of the composition rule
\begin{eqnarray}
l_{i_{1}}...\,l_{i_{m}}\!\blacktriangleleft\,l_{i_{m+1}}...\,l_{i_{m+n}}=\biguplus_{\substack{\sigma(1)<...<\sigma(m)\\\sigma(m+1)<...<\sigma(m+n)}}\{l_{i_{\sigma(1)}}...\,l_{i_{\sigma(m)}}l_{i_{\sigma(m+1)}}...\,l_{i_{\sigma(m+n)}}\}\ ,
\end{eqnarray}
where the index set runs over all permutations $\sigma$ of the set $\{1,2,...,m+n\}$ which preserve the relative order of $1,2,...,m$ and $m+1,m+2,...m+n$ respectively. One checks that a compatible decomposition rule is given by cutting a word into two parts and exchanging the prefix with the suffix, \textit{i.e.}
\begin{eqnarray}
\langle l_{i_1}...\,l_{i_k}\rangle=\biguplus_{j=0,...,k}\{(l_{i_{j+1}}...\,l_{i_{k}},l_{i_{1}}...\,l_{i_{j}})\}\ .
\end{eqnarray}
Note that this is the instance of a non-monoidal composition law.
Following the scheme of Section~\ref{Sect-Algebra} we arrive at the Hopf algebra structure given by the mappings
\begin{eqnarray}\label{Shuffle-mult}
l_{i_1}...\,l_{i_m}*\ l_{i_1}...\,l_{i_n}=\sum_{\substack{\sigma(1)<...<\sigma(m)\\\sigma(m+1)<...<\sigma(m+n)}}l_{i_{\sigma(1)}}...\,l_{i_{\sigma(m)}}l_{i_{\sigma(m+1)}}...\,l_{i_{\sigma(m+n)}}\ ,
\end{eqnarray}
\begin{eqnarray}\label{Shuffle-co-mult}
\Delta( l_{i_1}...\,l_{i_k})=\sum_{j=0}^{k}\ \ l_{i_{j+1}}...\,l_{i_{k}}\otimes\ l_{i_{1}}...\,l_{i_{j}}\ ,
\end{eqnarray}
\begin{eqnarray}
\varepsilon(l_{i_1}...\,l_{i_k})=0\ ,\ \ \ \ \ \ \ \ \varepsilon(\text{\O})=1\ ,
\end{eqnarray}
\begin{eqnarray}
S(l_{i_1}...\,l_{i_k})=(-1)^k\ l_{i_k}...\,l_{i_1}\ ,\ \ \ \ \ \ \ \ S(\text{\O})=\text{\O}\ .
\end{eqnarray}
We remark that for such constructed \emph{Shuffle algebra} the multiplication of Eq.~(\ref{Shuffle-mult}) is commutative and the co-product of Eq.~(\ref{Shuffle-co-mult}) is not co-commutative. 

\subsection{Graphs}

Let us consider a class of undirected \emph{graphs} $\mathcal{G}$ which one graphically represents as a collection of vertices connected by edges (we exclude isolated vertices). More formally a graph is defined as a mapping $\varGamma:E\longrightarrow V^{(2)}$ prescribing how the edges $E$ are attached to vertices $V$, where $V^{(2)}$ is a set of unordered pairs of vertices (not necessarily distinct); for a rigorous definition see \cite{Or90,Wi96,DiestelBook}. Let the size of a graph be the number of its edges $|\varGamma|=|E|$. 

An obvious composition rule in $\mathcal{G}$ consist in taking two graphs $\varGamma_2,\varGamma_1\in\mathcal{G}$ and drawing them one next to another, \textit{i.e.} 
\begin{eqnarray}
\varGamma_2\blacktriangleleft\varGamma_1=\varGamma_2\varGamma_1\ ,
\end{eqnarray}
where formally $\varGamma_2\varGamma_1:E_2\uplus E_1\longrightarrow V_2^{(2)}\uplus V_1^{(2)}$, such that $\varGamma_2\varGamma_1|_{E_2}=\varGamma_2$ and $\varGamma_2\varGamma_1|_{E_1}=\varGamma_1$.
\\
For example:\\
\begin{center}
\resizebox{0.7\columnwidth}{!}{\includegraphics{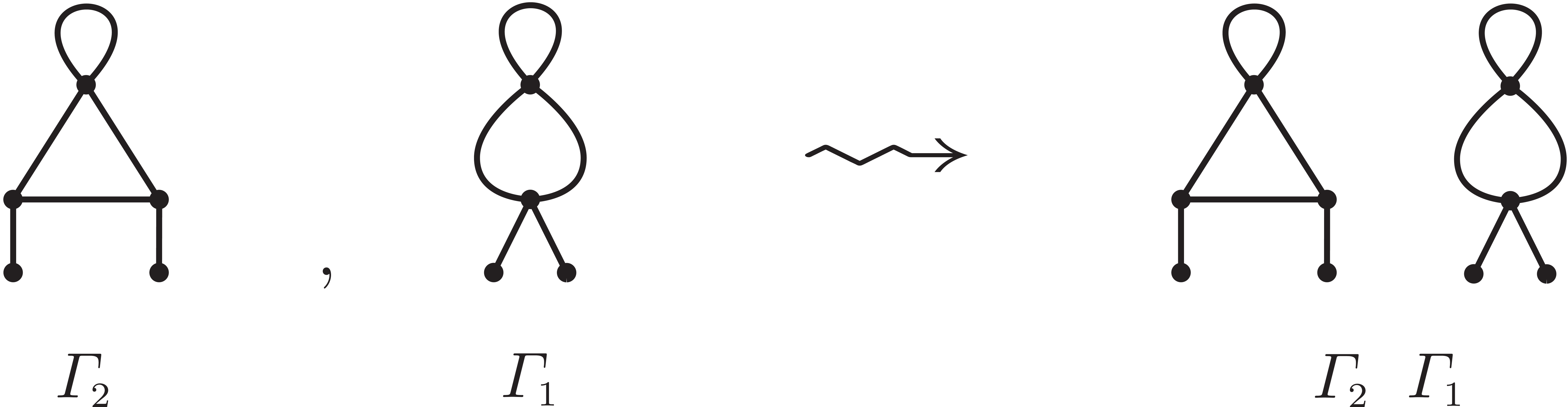}}
\end{center}
Such defined composition rule satisfies conditions \textbf{\texttt{(C1)}} -- \textbf{\texttt{(C4)}} and establishes a commutative algebra $\mathscr{G}$ (with void graph \O\ being the identity).

Observe that for a given graph $\varGamma\in\mathcal{G}$ each subset of its edges $L\subset E$ induces a subgraph $\varGamma|_L:L\longrightarrow V^{(2)}$ which is defined by restriction of $\varGamma$ to the subset $L$. Likewise, the remaining part of the edges $R=E-L$ gives rise to a subgraph $\varGamma|_R$. Thus, by considering ordered partitions of the set of edges into two subsets $L+R=E$, \textit{i.e.} $L\cup R=E$ and $L\cap R=\emptyset$, we end up with pairs $(\varGamma|_L,\varGamma|_R)$ of disjoint graphs. 
For example:\\
\begin{center}
\resizebox{0.6\columnwidth}{!}{\includegraphics{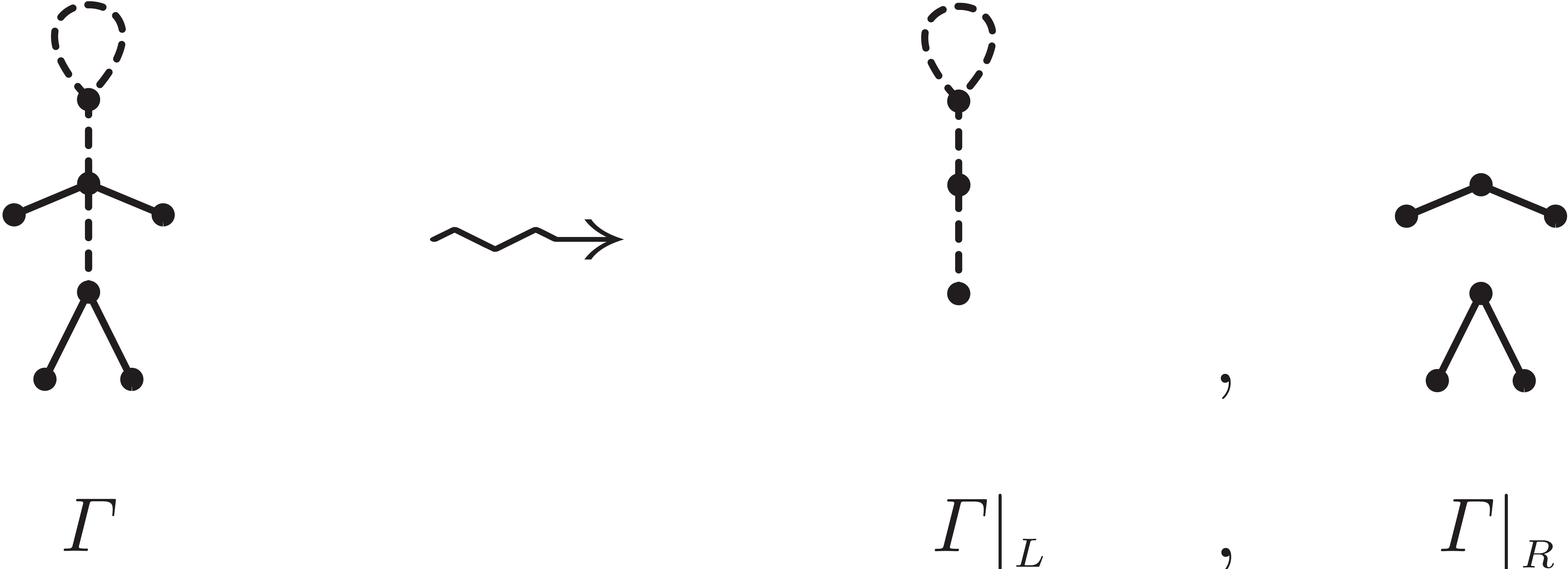}}
\end{center}
is a decompositions of a graph $\varGamma\leadsto(\varGamma|_L,\varGamma|_R)$ induced by the choice of the dashed edges $L\subset E$ which make the left component in the splitting.
This suggests the following definition of the decomposition rule
\begin{eqnarray}
\langle\varGamma\rangle=\biguplus_{L+R=E}\{(\varGamma|_L,\varGamma|_R)\}\ .
\end{eqnarray}
One checks that conditions \textbf{\texttt{(D1)}} -- \textbf{\texttt{(D5)}} and \textbf{\texttt{(CD2)}} hold, and we obtain a graded Hopf algebra $\mathscr{G}$ with the grading given by the number of edges. Its structure is given by
\begin{eqnarray}
\varGamma_2*\varGamma_1=\varGamma_2\varGamma_1\ ,
\end{eqnarray}
\begin{eqnarray}
\Delta(\varGamma)=\sum_{L+R=E}\varGamma|_L\otimes\varGamma|_R\ ,
\end{eqnarray}
\begin{eqnarray}
\varepsilon(\varGamma)=0\ ,\ \ \ \ \ \ \ \ \varepsilon(\text{\O})=1\ ,
\end{eqnarray}
\begin{eqnarray}
S(\varGamma)=\sum_{\substack{L_1+...+L_n=E\\L_1,...,L_n\neq\,\text{\O}}}(-1)^n\ \ \varGamma|_{L_1}...\,\varGamma|_{L_n}\ , \ \ \ \ \ \ \ \ \ \ \ S(\text{\O})=\text{\O}\ .
\end{eqnarray}
So defined \emph{algebra of graphs} is both commutative and co-commutative.







\subsection{Trees and Forests}
A \emph{rooted tree} is a graph without cycles with one distinguished vertex, called the root. Let $\mathcal{T}$ denote the class of rooted trees. A \emph{forest} is a collection of rooted trees and the pertaining combinatorial class has the specification $\mathcal{F}=\textsc{MSet}(\mathcal{T})$. Size of a tree (forest) is defined as the number of vertices. 

We will consider class $\mathcal{F}$ and define composition of forests as a multiset union (like for graphs), \textit{i.e.}
\begin{eqnarray}
\varGamma_2\blacktriangleleft\varGamma_1=\varGamma_2\varGamma_1\ .
\end{eqnarray}
For example:\\
\begin{center}
\resizebox{0.7\columnwidth}{!}{\includegraphics{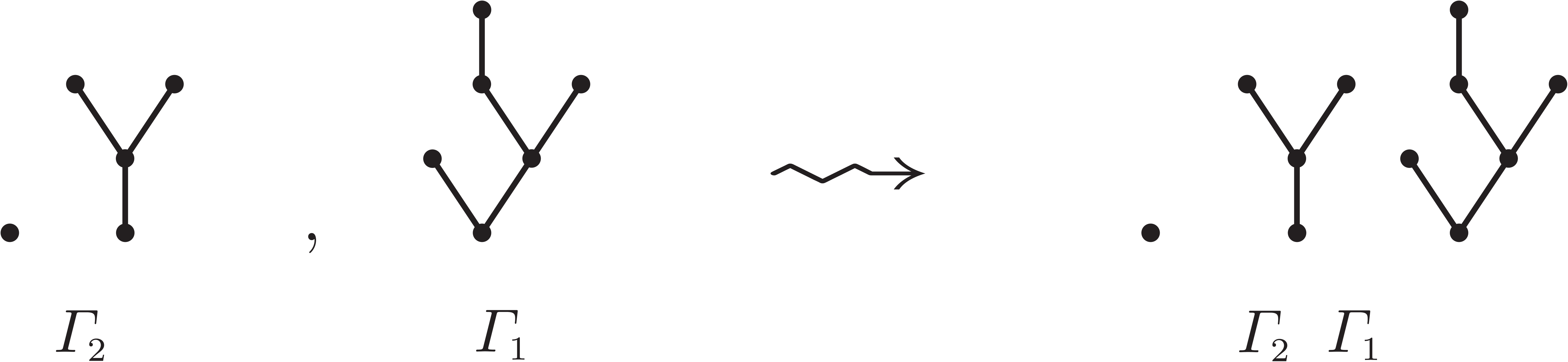}}
\end{center}
Note that $(\mathcal{F},\blacktriangleright)$ is a (commutative) monoid generated by the rooted trees $\mathcal{T}$.

For a given a tree $\tau\in\mathcal{T}$ one distinguishes subtrees $\tau^{r}\subset\tau$ which share the same root with $\tau$, called proper subtrees (the empty tree \O\ is considered as a proper subtree as well). Observe that the latter obtains by trimming $\tau$ to the required shape $\tau^{r}$, and the branches which are cut off form a forest of trees denoted by $\tau^{c}$ (with the roots next to the cutting). Decomposition of a tree is defined as any splitting $\tau\leadsto(\tau^{c},\tau^{r})$ into a pair consisting of a proper subtree taken in the second component and the remaining forest in the first one. In other words
\begin{eqnarray}\label{tree-decomp}
\langle\tau\rangle=\biguplus_{\tau^r\subset\tau}\{(\tau^c,\tau^r)\}\ ,
\end{eqnarray}
where the disjoint union ranges over proper subtrees $\tau^r$ of $\tau$, and $\tau^c$ is a forest of trees which 'complements' $\tau^r$ to $\tau$. For example:
\\
\begin{center}
\resizebox{1\columnwidth}{!}{\includegraphics{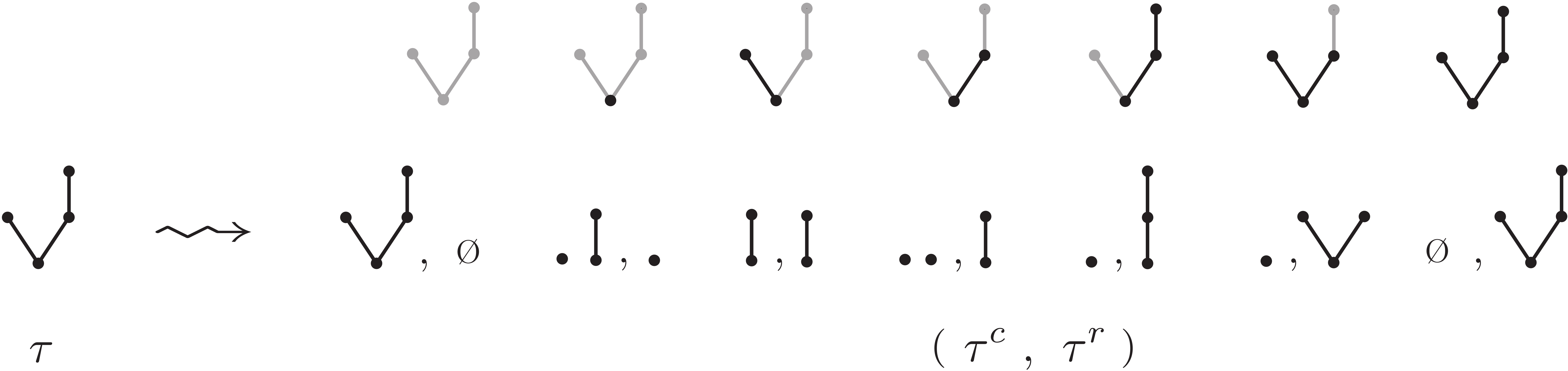}}
\end{center}
enumerates all possible decompositions of a tree (at the top row proper subtrees $\tau^r$ are drawn in black while the completing forests $\tau^c$ are drawn in gray).
\\
Since trees $\mathcal{T}$ generate forests $\mathcal{F}$, we extend the decomposition rule to any forest $\varGamma=\tau_n\ ...\ \tau_1\in\mathcal{F}$ using Eq.~(\ref{extension-decomp}) and obtain
\begin{eqnarray}\label{forest-decomp}
\langle\tau_n\ ...\ \tau_1\rangle=\biguplus_{\tau_n^r\subset\tau_n\,,\,...\,,\,\tau_1^r\subset\tau_1}\{(\tau_n^c\ ...\ \tau_1^c,\tau_n^r\ ...\ \tau_1^r)\}\ ,
\end{eqnarray}
which comes down to trimming some of the branches off the whole forest and gathering them in the first component $\varGamma^c=\tau_n^c\ ...\ \tau_1^c$ whilst keeping the rooted parts in the second one $\varGamma^r=\tau_n^r\ ...\ \tau_1^r$. Hence, we will briefly write
\begin{eqnarray}\label{forest-decomp}
\langle\varGamma\rangle=\biguplus_{\varGamma^r\subset\varGamma}\{(\varGamma^c,\varGamma^r)\}\ .
\end{eqnarray}
Following the construction of Section~\ref{Sect-Algebra} one obtains a graded Hopf algebra $\mathscr{F}$ with the grading given by the number of vertices.
The required mappings take the form
\begin{eqnarray}
\varGamma_2*\varGamma_1=\varGamma_2\ \varGamma_2\ ,
\end{eqnarray}
\begin{eqnarray}
\Delta(\varGamma)=\sum_{\varGamma^r\subset\varGamma}\varGamma^c\otimes\varGamma^r\ ,
\end{eqnarray}
\begin{eqnarray}
\varepsilon(\varGamma)=0\ ,\ \ \ \ \ \ \ \ \varepsilon(\text{\O})=1\ ,
\end{eqnarray}
\begin{eqnarray}
S(\varGamma)=-\sum_{\varGamma^r\subset\varGamma}S(\varGamma^c)\ \varGamma^r\ ,\ \ \ \ \ \ \ \ S(\text{\O})=\text{\O}\ .
\end{eqnarray}
Such constructed \emph{algebra of forests} $\mathscr{F}$ is commutative but not co-commutative. 
We remark that this Hopf algebra was first introduced by J. C. Butcher \cite{Butcher1972,Brouder2004} and recently it was rediscovered by A. Connes and D. Kreimer \cite{Kreimer1998,ConnesKreimer1998} in the context of renormalization in quantum field theory.

\acknowledgements

I wish to thank Gerard Duchamp, Philippe Flajolet, Andrzej Horzela, Karol A. Penson and Allan I. Solomon for important discussions on the subject.
Most of this research was carried out in the Mathematisches Forschungsinstitut Oberwolfach (Germany) and the Laboratoire d'Informatique de l'Universit\'e Paris--Nord in Villetaneuse (France) whose warm hospitality is greatly appreciated.
The author acknowledges support from the Agence Nationale de la Recherche under the programme no. ANR-08-BLAN-0243-2 and the Polish Ministry of Science and Higher Education grant no. N202 061434.

\bibliographystyle{alpha}
\bibliography{Bibliography}

\newcommand{\etalchar}[1]{$^{#1}$}
\begin{thebibliography}{BDH{\etalchar{+}}10}

\bibitem[Abe80]{Ab80}
E.~Abe.
\newblock {\em Hopf Algebras}.
\newblock Cambridge University Press, 1980.

\bibitem[BDH{\etalchar{+}}10]{BlDuHoPeSo10}
P.~Blasiak, G.~H.~E. Duchamp, A.~Horzela, K.~A. Penson, and A.~I. Solomon.
\newblock Combinatorial {A}lgebra for second-quantized {Q}uantum {T}heory.
\newblock {\em Adv. Theor. Math. Phys.}, 14(4), 2010.
\newblock Article in Press, arXiv:1001.4964 [math-ph].

\bibitem[BLL98]{BergeronBook}
F.~Bergeron, G.~Labelle, and P.~Leroux.
\newblock {\em Combinatorial Species and Tree-like Structures}.
\newblock Cambridge University Press, 1998.

\bibitem[Bou89]{Bo89}
N.~Bourbaki.
\newblock {\em Algebra I: Chapters 1--3}.
\newblock Springer-Verlag, 1989.

\bibitem[Bro04]{Brouder2004}
C.~Brouder.
\newblock Trees, {R}enormalization and {D}ifferential {E}quations differential
  equations.
\newblock {\em BIT Num. Math.}, 44:425--438, 2004.

\bibitem[But72]{Butcher1972}
J.~C. Butcher.
\newblock An {A}lgebraic {T}heory of {I}ntegration {M}ethods.
\newblock {\em Math. Comput.}, 26:79--106, 1872.

\bibitem[Car07]{CartierHopfPrimer}
P.~Cartier.
\newblock A {P}rimer of {H}opf {A}lgebrasrimer of hopf algebras.
\newblock In {\em Frontiers in Number Theory, Physics, and Geometry II}, pages
  537--615. Springer, 2007.

\bibitem[CK98]{ConnesKreimer1998}
A.~Connes and D.~Kreimer.
\newblock Hopf {A}lgebras, {R}enormalization and {N}oncommutative {G}eometry.
\newblock {\em Commun. Math. Phys.}, 199:203--242, 1998.

\bibitem[Die05]{DiestelBook}
R.~Diestel.
\newblock {\em Graph Theory}.
\newblock Springer-Verlag, 3rd edition, 2005.

\bibitem[FS09]{FlajoletBook}
P.~Flajolet and R.~Sedgewick.
\newblock {\em Analytic Combinatorics}.
\newblock Cambridge University Press, 2009.

\bibitem[GJ83]{GouldenBook}
I.~P. Goulden and D.~M. Jackson.
\newblock {\em Combinatorial Enumeration}.
\newblock John Wiley \& Sons, 1983.

\bibitem[GL89]{GrLa89}
R.~Grossman and R.~G. Larson.
\newblock Hopf-{A}lgebraic {S}tructure of {F}amilies of {T}rees.
\newblock {\em J. Algebra}, 126(1):184--210, 1989.

\bibitem[Joy81]{Jo81}
A.~Joyal.
\newblock Une th{\'e}orie combinatoire des s{\'e}ries formelles.
\newblock {\em Adv. Math.}, 42(1):1--82, 1981.

\bibitem[JR79]{JoniRota}
S.~A. Joni and G.~C. Rota.
\newblock Coalgebras and {B}ialgebras in {C}ombinatorics.
\newblock {\em Stud. Appl. Math.}, 61:93--139, 1979.

\bibitem[Kre98]{Kreimer1998}
D.~Kreimer.
\newblock On the {H}opf algebra structure of perturbative quantum field theory.
\newblock {\em Adv. Theor. Math. Phys.}, 2:303--334, 1998.

\bibitem[Lot83]{LothaireBook}
M.~Lothaire.
\newblock {\em Combinatorics on Words}.
\newblock Addison-Wesley, 1983.

\bibitem[Ore90]{Or90}
O.~Ore.
\newblock {\em Graphs and their Uses}.
\newblock Mathematical Association of America, 2nd edition, 1990.

\bibitem[Reu93]{ReutenauerBook}
C.~Reutenauer.
\newblock {\em Free Lie Algebras}.
\newblock Oxford University Press, 1993.

\bibitem[Swe69]{SweedlerBook}
M.~E. Sweedler.
\newblock {\em Hopf Algebras}.
\newblock Benjamin, 1969.

\bibitem[Wil96]{Wi96}
R.~J. Wilson.
\newblock {\em Introduction to Graph Theory}.
\newblock Addison-Wesley, 4th edition, 1996.

\end{thebibliography}

\end{document}